\documentclass[11pt]{amsart}
\usepackage[latin1]{inputenc}
\usepackage{hyperref}
\usepackage{amsfonts}
\usepackage{amssymb}
\usepackage{amsmath}
\usepackage{amsthm}
\usepackage{latexsym}
\usepackage{color}
\usepackage{amscd}
\usepackage{graphicx}
\usepackage{enumitem} 
\usepackage[margin=1.33in]{geometry}
\usepackage{esint} 
\usepackage{epsfig}

\allowdisplaybreaks

\numberwithin{equation}{section}
\newcounter{exercise}

\newtheorem{theorem}{Theorem}

\theoremstyle{definition}

\theoremstyle{remark}

\newcommand\R{{\mathbb R}}
\newcommand\T{{\mathbb T}}

\newcommand{\cC}{\mathcal C}

\newcommand{\cQ}{\mathcal Q}

\newcommand{\dd}{{\, \mathrm d}}

\let\oldmarginpar\marginpar
\renewcommand\marginpar[1]{\-\oldmarginpar[\raggedleft\footnotesize #1]%
{\raggedright\footnotesize #1}}

\def\signcm{\bigskip \begin{center} {\sc Cl\'ement
      Mouhot\par\vspace{3mm}
      University of Cambridge\par
      DPMMS, Centre for Mathematical Sciences\par
      Wilberforce road, Cambridge CB3 0WA, UK
      \par\vspace{3mm} e-mail:}
    \tt{C.Mouhot@dpmms.cam.ac.uk} \end{center}}

\begin{document}

\title[De Giorgi-Nash-Moser and H\"ormander theories: new
interplay]{De Giorgi-Nash-Moser and H\"ormander theories: new
  interplay}

\author{C. Mouhot}

\date{\today}

\begin{abstract}
  We report on recent results and a new line of research at the
  crossroad of two major theories in the analysis of partial
  differential equations. The celebrated De~Giorgi-Nash-Moser theory
  shows H\"older estimates and the Harnack inequality for uniformly
  elliptic or parabolic equations with rough coefficients in
  divergence form. The theory of hypoellipticity of H\"ormander shows,
  under general ``bracket'' conditions, the regularity of solutions to
  partial differential equations combining first and second order
  derivative operators when ellipticity fails in some directions. We
  discuss recent extensions of the De Giorgi-Nash-Moser theory to
  hypoelliptic equations of Kolmogorov (kinetic) type with rough
  coefficients. These equations combine a first-order skew-symmetric
  operator with a second-order elliptic operator involving derivatives
  in only certain variables, and with rough coefficients. We then
  discuss applications to the Boltzmann and Landau equations in
  kinetic theory and present a program of research with some open
  questions.
\end{abstract}
\maketitle

\tableofcontents

\section{Introduction}
\label{sec:intro}

\subsection{Kinetic theory}

Modern physics goes back to Newton and classical mechanics, and was
later expanded into the understanding of electric and magnetic forces
were (Amp\`ere, Faraday, Maxwell), large velocities and large scales
(Lorentz, Poincar\'e, Minkowski, Einstein), small-scale particle
physics and quantum mechanics (Planck, Einstein, Bohr, Heisenberg,
Born, Jordan, Pauli, Fermi, Schr\"odinger, Dirac, De Broglie, Bose,
etc.). However, all these theories are classically devised to study
one physical system (planet, ship, motor, battery, electron,
spaceship, etc.) or a small number of systems (planets in the Solar
system, electrons in a molecule, etc.)  In many situations though, one
needs to deal with an assembly made up of elements so numerous that
their individual tracking is not possible: galaxies made of hundreds
of billions of stars, fluids made of more than $10^{20}$ molecules,
crowds made of thousands of individuals, etc. Taking such large
numbers into account leads to new effective laws of physics, requiring
different models and concepts.  This passage from microscopic rules to
macroscopic laws is the founding principle of statistical physics. All
branches of physics (classical, quantum, relativistic, etc.)  can be
studied from the point of view of statistical physics, in both
stationary and dynamical perspectives. It was first done with the laws
of classical mechanics, which resulted in kinetic theory, discovered
by Maxwell \cite{Maxwell1867} and Boltzmann \cite{Boltzmann1872} in
the 19th century after precursory works by D.~Bernoulli, Herapath,
Waterston, Joule, K\"onig, Clausius.

Kinetic theory replaces a huge number of objects, whose physical
states are described by a certain phase space, and whose properties
are otherwise identical, by a {\em statistical distribution} over that
phase space.  The fundamental role played by the velocity (kinetic)
variable inaccessible to observation was counter-intuitive, and
accounts for the denomination of {\em kinetic} theory. The theory
introduces a distinction between three scales: the macroscopic scale
of phenomena which are accessible to observation; the microscopic
scale of molecules and infinitesimal constituents; and an intermediate
scale, loosely defined and often called {\em mesoscopic}. This is the
scale of phenomena which are not accessible to macroscopic observation
but already involve a large number of particles, so that statistical
effects are significant.

\subsection{Main equations of kinetic theory}

Maxwell wrote the first (weak) form of the evolution equation known
now as the {\bf Boltzmann equation}: the unknown is a (non-negative)
density function $f(t,x,v)$, standing for the density of particles at
time $t$ in the phase space $(x,v)$; the equation, in modern writing
and assuming the absence of external forces, is
\begin{equation}\label{MBE}
  \frac{\partial f}{\partial t} + v\cdot\nabla_x f
  = Q(f,f).
\end{equation}
The left-hand side describes the evolution of $f$ under the action of
transport, with the \emph{free streaming operator}. The right-hand
side describes elastic collisions, with the nonlinear
\textbf{Boltzmann collision operator}:
\begin{equation}\label{BCO}
  Q(f,f) =
  \int_{\mathbb{R}^{3}}\,\int_{\mathbb S^{2}} B(v-v_*,\omega)
  \Bigl(f(t,x,v^{\prime })\,f(t,x,v_{\ast }^{\prime})
  -f(t,x,v)\,f(t,x,v_{\ast })\Bigr) \dd v_{\ast } \dd \omega.
\end{equation}
Note that this operator is localized in $t$ and $x$, quadratic, and
has the structure of a tensor product with respect to
$f(t,x,\cdot)$. The velocities $v'$ and $v'_*$ should be thought of as
the velocities of a pair of particles before collision, while $v$ and
$v_*$ are the velocities after that collision: the formulas are
$v' = v- \langle v-v_*,\omega\rangle \,\omega$ and
$v'_* = v_* + \langle v-v_*,\omega\rangle\,\omega$. When one computes
$(v,v_*)$ from $(v',v'_*)$ (or the reverse), conservation laws of the
mass, momentum and energy are not enough to yield the result, with
only 4 scalar conservation laws for 6 degrees of freedom. The unit
vector $\omega \in \mathbb S^2$ removes this ambiguity: in the case of
colliding hard spheres, it can be thought of as the direction of the
line joining the two centers of the particles.  The kernel
$B(v-v_*,\omega)$ describes the relative frequency of vectors
$\omega$, depending on the relative impact velocity $v-v_*$; it only
depends on the modulus $|v-v_*|$ and the deflection angle $\theta$
between $v-v_*$ and $v'-v'_*$. Maxwell computed it for hard spheres,
where $B\sim |v-v_*|\sin\theta$, and for inverse power forces, where
the kernel factorizes as $B \sim |v-v_*|^\gamma b(\cos \theta)$.
Maxwell showed that if the force is repulsive, proportional to
$r^{-\alpha}$ ($r$ the inter-particle distance), then
$\gamma = (\alpha-5)/(\alpha-1)$ and
$b(\cos \theta)\simeq \theta^{-(1+2s)}$ as $\theta\to 0$, where
$2s = 2/(\alpha-1)$. In particular, the kernel is usually {\em
  non-integrable} as a function of the angular variable: this is a
general feature of \emph{long-range interactions}, nowadays sometimes
called ``non-cutoff property''.

The case $\alpha=5$, $\gamma=0$ and $2s=1/2$ is called \emph{Maxwell
  molecules} \cite{Maxwell1867}, the case $\alpha \in (5,+\infty)$,
$\gamma >0$ and $2s \in (0,1/2)$ is called \emph{hard potentials
  (without cutoff)}, the case $\alpha \in [3,5)$, $\gamma \in [-1,0)$,
$2s \in (1/2,1]$ is called \emph{moderately soft potentials (without
  cutoff)}, and finally the case $\alpha \in (2,3)$,
$\gamma \in (-3,-1)$, $2s \in (1,2)$ is called \emph{very soft
  potentials (without cutoff)}. The limits between hard and soft
potentials ($\gamma =0$) and between moderately and very soft
potentials ($\gamma + 2s =0$) are commonly taken as defining the
``hard'' / ``moderately soft'' / ``very soft'' terminology in any
dimension, for kernel of the form $B = |v-v_*|^\gamma b(\cos \theta)$
with $b(\cos \theta)\simeq \theta^{-(1+2s)}$.

In order to find the stationary solutions, that is, time-independent
solutions of \eqref{BCO}, the first step is to identify particular
\emph{hydrodynamic} density functions, which make the collision
contribution vanish: these are {\em Gaussian distributions with a
  scalar co-variance}
$f(v) = \rho\,(2\pi T)^{-3/2} \, e^{-\frac{|v-u|^2}{2T}}$, where the
parameters $\rho>0$, $u\in\R^3$ and $T>0$ are the local density, mean
velocity, and temperature of the fluid. These parameters can be fixed
throughout the whole domain (providing in this case an
\emph{equilibrium distribution}), or depend on the position $x$ and
time $t$; in both cases collisions will have no effect. As pointed out
in Maxwell's seminar paper, and later proved rigorously at least in
some settings
\cite{MR1115587,MR1213991,MR2025302,MR2197021,MR2517786}, the
Boltzmann equation is connected to classical fluid mechanical
equations on $\rho$, $u$ and $T$, and leads to them in certain
regimes. This provides a rigorous connection between the mesoscopic
(kinetic) level and the macroscopic level. At the other end of the
scales, a rigorous derivation of the Boltzmann equation from many-body
Newtonian mechanics for short time and short-range interactions was
obtained by Lanford \cite{MR0479206} for hard spheres; see also
\cite{MR2625983} for an extension to more general short-range
interactions, and the recent works \cite{MR3157048,MR3190204} that
revisit and complete the initial arguments of Lanford and King. Note
however that at the moment the equivalent of Lanford theorem for the
Boltzmann equation with long-range interactions is still missing, see
\cite{MR3607474} for partial progresses.

To summarise the key mathematical points: the Boltzmann equation is an
integro-(partial)-differential equation with non-local operator in the
kinetic variable $v$. Moreover for long-range interactions with
repulsive force $F(r) \sim r^{-\alpha}$, this non-local operator has a
singular kernel and shows, as we will see, fractional ellipticity of
order $2/(\alpha -1)$. The Boltzmann equation ``contains'' the
hydrodynamic, and it is a \emph{fundamental} equation in the sense
that it is derived rigorously, at least in some settings, from
microscopic first principles. From now on, we consider the position
variable in $\R^3$ or in the periodic box $\T^3$.

In the limit case $s \to 1$ (the Coulomb interactions), the Boltzmann
collision operator is ill-defined. Landau \cite{Landau1936} proposed
an alternative operator for these Coulomb interactions that is now called
the \textbf{Landau-Coulomb operator}
\begin{align*}
  Q(f,f) = \nabla_v \cdot \left( \int_{\R^3}
  {\bf P}_{(v-v_*)^\bot} \Big( f(t,x,v_*) \nabla_v f(t,x,v) - f(t,x,v)
  \nabla_v f(t,x,v_*) \Big) |v-v_*|^{\gamma+2} \dd v_* \right) 
\end{align*}
where ${\bf P}_{(v-v_*)^\bot}$ is the orthogonal projection along
$(v-v_*)^\bot$ and $\gamma = -3$. It writes also
\begin{align}
  \label{eq:LCop} 
  Q(f,f) = \nabla_v \cdot \big( A[f] \nabla_v
  f + B[f] f \big) 
\end{align}  
\begin{equation*}
\mbox{with }
\begin{cases}
  \displaystyle A[f] (v) = \int_{\mathbb R^3} \left( I -
    \frac{w}{|w|}\otimes \frac{w}{|w|}\right) \,
  |w|^{\gamma+2} \, f(t,x,v-w) \dd w,\\[3mm]
  \displaystyle B[f](v) = - \int_{\mathbb R^3} |w|^{\gamma} \, w \,
  f(t,x,v-w) \dd w.
\end{cases}
\end{equation*}
This operator is a nonlinear drift-diffusion operator with
coefficients given by convolution-like averages of the unknown. This
is a non-local integro-differential operator, with second-order local
ellipticity. The resulting {\bf Landau
  equation}~\eqref{MBE}-\eqref{eq:LCop} again ``contains'' the
hydrodynamic. It is also considered \emph{fundamental} because of its
closed link to the Boltzmann equation for Coulomb interactions (note
however that the equivalent to Lanford theorem for the Landau equation
is lacking, even at a formal level, see \cite{MR3040372} for partial
progresses). Because of the difficulty to handle the very singular
kernel of the Landau-Coulomb operator, it is common to introduce
artificially a scale of models by letting $\gamma$ vary in $[-3,1]$
(or even $[-d,1]$ in general dimension $d$). The terminology
\emph{hard potentials}, \emph{Maxwell molecules}, \emph{soft
  potentials} are used as for the Boltzmann collision operator when
$\gamma >0$, $\gamma=0$, $\gamma <0$ respectively. The terminology
\emph{moderately soft potentials} corresponds here (since $s=1$) to
$\gamma \in (-2,0)$.

\subsection{Open problems and conjectures}

\subsubsection{The Cauchy problem}

The first mathematical question when studying the previous fundamental
kinetic equations (Boltzmann and Landau equations) is the Cauchy
problem, i.e. existence, uniqueness and regularity of
solutions. Short-time solutions have been constructed, as well as
global solutions close to the trivial stationary solution or with
space homogeneity: see \cite{GMM} and the references therein for some
of the most recent results for the Boltzmann equation with short-range
interactions, see
\cite{MR2863853,MR2793203,MR2679369,MR2795331,MR2784329} for the
Boltzmann equation with long-range interactions, and see
\cite{MR1946444} for the Landau equation. However the construction of
solutions ``in the large'' remains a formidable open problem. Since
weak ``renormalised'' solutions have been constructed by DiPerna and
Lions \cite{MR1014927} that play a similar role to the Leray
\cite{MR1555394} solutions in fluid mechanics, this open problem can
be compared with the millennium problem of the regularity of solutions
to 3D incompressible Navier-Stokes equations.

\subsubsection{Study of a priori solutions}

Given that the Cauchy problem in the large seems out of reach at the
moment, Truesdell and Muncaster \cite{MR554086} remarked almost 40
years ago that: ``\textit{Much effort has been spent toward proof that
  place-dependent solutions exist for all time. [\dots] The main
  problem is really to discover and specify the circumstances that
  give rise to solutions which persist forever. Only after having done
  that can we expect to construct proofs that such solutions exist,
  are unique, and are regular.}'' In other words, the $H$-theorem and
the mathematical understanding of irreversibility are so important in
the theory of Maxwell and Boltzmann that it cannot wait for the
tremendously difficult issue of global well-posedness to be
settled. Cercignani then formulated a precise conjecture along this
idea, postulating in \cite{MR715658} a linear relation between the
entropy production functional and the relative entropy functional of
any a priori given classical solutions. The resolution of this
conjecture, for certain interactions, lead to precise new quantitative
information on a priori solutions of the Boltzmann and Landau equation
(see \cite{MR2116276,MR2765747,MR2197542,GMM,MR3625186}). And it lead
to the related question of the optimal relaxation rates of a priori
solutions, with minimal regularity and moments conditions. It is now
fairly well understood for many interactions. The results obtained
along this line of research can all be summarised into the following
general form: \medskip

\noindent {\bf Conditional relaxation.} Any solution to the Boltzmann
(resp. Landau) equation in
$L^\infty_x (\T^3; L^1_v(\R^3, (1+|v|)^k \dd v))$, $k >2$ (or a
closely related functional space as large as possible) converges to
the thermodynamic equilibrium with the optimal rate dictated by the
linearized equation.  \medskip

Note however that an interesting remaining open question in this
program is to obtain a result equivalent to \cite{GMM,MR3625186} in
the case of the Boltzmann equation with long-range interactions (with
fractional ellipticity in the velocity variable).

\subsubsection{Regularity conjectures for long-range interactions}

In the case of long-range interactions, the Boltzmann and
Landau-Coulomb operators show local ellipticity provided the solution
enjoys some pointwise bounds on the hydrodynamic fields
$\rho(t,x):= \int_{\R^3} f \dd v$,
$e(t,x):= \int_{\R^3} f |v|^2 \dd v$ and the local entropy
$h(t,x):= \int_{\R^3} f \ln f \dd v$. Whereas it is clear in the case
of the Landau-Coulomb operator, it was understood almost two decades
ago in the case of the Boltzmann collision operator
\cite{MR1649477,MR1715411,MR1765272}. This had lead colleagues working
on non-local operators and fully nonlinear elliptic problems like
Silvestre and Guillen and co-authors to attempt to use barriers'
techniques reminiscent to the Krylov-Safonov theory \cite{MR563790} in
order to obtain pointwise bounds for solutions to these
equations. These first attempts, while unsuccessful, later proved
crucial in attracting the attention of a larger community on this
problem. And these authors rapidly reformulated the initial goal into,
again, \emph{conditional} conjectures on the regularity of the form:
\medskip

\noindent {\bf Conditional regularity.} Consider any solution to the
Boltzmann equation with long-range interactions (resp. Landau
equation) on a time interval $[0,T]$ such that its hydrodynamic
fields are bounded:
\begin{equation}\label{eq:hydro}
  \forall \, t \in [0,T], \ x \in \T^3, \quad m_0 \le \rho(t,x) \le m_1, \quad e(t,x) \le e_1, \quad  h(t,x) \le h_1
\end{equation}
where $m_0$, $m_1$, $e_1$, $h_1 >0$. Then the solution is bounded and
smooth on $(0,T]$.  \medskip

Note that this conjecture can be strengthened by removing the
assumption that the mass is bounded from below and replacing it by a
bound from below on the total mass
$\int_{\T^3} \rho(t,x) \dd x \ge M_0>0$. Mixing in velocity through
collisions combined with transport effects indeed generate lower
bounds in many settings, see
\cite{MR2153518,MR2746671,MR3375544,MR3356579}; moreover it was indeed
proved for the Landau equation with moderately soft potentials in
\cite{henderson2017local}.

This conjecture is now been partially solved in the case of the Landau
equation, when the interaction is ``moderately soft''
$\gamma \in (-2,0)$. This result has been the joint efforts of several
groups \cite{gimv,henderson2017c,henderson2017local,IM-toy}, and this is the
object of the next section. It is currently an ongoing program of
research in the case of the Boltzmann equation with hard and
moderately soft potentials, and this is the object of the fourth and
last section. The conjecture interestingly remains open in the case of
very soft potentials for both equations, and making progress in this
setting is likely to require new conceptual tools.

\section{De Giorgi-Nash-Moser meet H\"ormander}
\label{sec:dgnm}

\subsection{The resolution of Hilbert $19$-th problem} The De
Giorgi-Nash-Moser theory \cite{DeG56,DeG,nash,MR0170091,moser} was
born out of the attempts to answer Hilbert's 19th problem. This
problem is about proving the analytic regularity of the minimizers $u$
of an energy functional $\displaystyle \int_U L(\nabla u) \dd x$, with
$u : \R^d \to \R$ and where the \emph{Lagrangian}
$L: \mathbb R^d \to \mathbb R$ satisfies growth, smoothness and
convexity conditions and $U \subset \R^d$ is a compact domain. The
Euler-Lagrange equations for the minimizers take the form
\[
  \nabla \cdot \Big[ \nabla L (\nabla u) \Big] = 0 \quad \mbox{i.e.}
  \quad \sum_{i,j=1} ^d \underbrace{\big[(\partial_{ij} L)(\nabla
    u)\big]}_{b_{ij}} \partial_{ij}u = \sum_{i,j=1} ^d b_{ij} \partial_{ij}
  u = 0.
\]
For instance the Dirichlet energy $L(p) = |p|^2$ leads to linear
Euler-Lagrange equations, whereas the minimal surface energy
$L(p) = \sqrt{1+|p|^2}$ leads to nonlinear Euler-Lagrange
equations. With suitable assumptions on $L$ and the domain, the
pointwise control of $\nabla u$ was known in the 1950s. However
applying the Schauder estimates to get higher regularity requires more
information: if $u \in C^{1,\alpha}$ with $\alpha >0$ then
$b_{ij} \in C^\alpha$ and Schauder estimates \cite{MR1545448} imply
$u \in C^{2,\alpha}$; a bootstrap argument then yields higher
regularity, and analyticity follows from this $C^\infty$ regularity
\cite{MR1511259,MR0001425}.

Hence, apart from specific result in two dimensions \cite{MR1501936},
the missing piece in solving Hilbert 19th problem, in the 1950s, was
the proof of the H\"older regularity of $\nabla u$. The equation
satisfied by a derivative $f := \partial_k u$ is the divergence form
elliptic equation:
\[
  \sum_{i,j=1}^d \partial_i \Big[ \underbrace{(\partial_{ij} L)(\nabla
    u)}_{a_{ij}} \partial_j f \Big] =\nabla \cdot \left( A \nabla f
  \right) = 0.
\]
De Giorgi \cite{DeG} and Nash \cite{nash} independently proved this
H\"older regularity of $f$ under the sole assumption that the
symmetric matrix $A:=(a_{ij})$ satisfies the controls
$0 < \lambda \le A \le \Lambda$, and is measurable (no regularity is
assumed). The proof of Nash uses what is now called the ``Nash
inequality'', an $L \log L$ energy estimate, and refined estimates on
the fundamental solution. The proof of De Giorgi uses an iterative
argument to gain integrability, and an ``isoperimetric argument'' to
control how oscillations decays when refining the scale of
observation. Moser later gave an alternative proof
\cite{MR0170091,moser} based on one hand on an iterative gain of
integrability, formulated differently but similar to that of De
Giorgi, and on the other hand on relating Lebesgue norms on $f$ and
$1/f$ through energy estimates on the equation satisfied by $g:=\ln f$
and the use of a \emph{Poincar\'e inequality}; the proof of Moser had
an important further contribution in that it also proved, as an
intermediate step towards the Hölder regularity, the \emph{Harnack
  inequality} for the equation considered, i.e. a universal control on
the ratio between local maxima and local minima.

Let us mention that the De Giorgi-Nash-Moser (DGNM) theory only
considers elliptic or parabolic equations in \emph{divergence
  form}. An important counterpart result for non-divergence elliptic
and parabolic equations was later discovered by Krylov and Safonov
\cite{MR563790}. The extension of the DGNM theory to hypoelliptic
equations with rough coefficients that we present in this section
requires the equation to be in divergent form. It is an open problem
whether the Krylov-Safonov theory extends to hypoelliptic
non-divergent equations of the form discussed below.

\subsection{The theory of hypoellipticity} The DGNM theory has
revolutionised the study of nonlinear elliptic and parabolic partial
differential equations (PDEs). However it remained limited to PDEs
where the diffusion acts in all directions of the phase space. In
kinetic theory, as soon as the solution is non spatially homogeneous,
the diffusion or fractional diffusion in velocity is combined to a
conservative Hamiltonian dynamic in position and velocity. This
structure is called \textbf{hypoelliptic}.

The study of regularity properties of such equations can be traced
back, at the linear level, to the short note of Kolmogorov
\cite{MR1503147}. This note considered the combination of free
transport with drift-diffusion in velocity: the law satisfies what is
now sometimes called the \emph{Kolmogorov equation}, that writes
$\partial_t f + v \cdot \partial_x f = \Delta_v f$ on
$x,v \in \R^d$ in the simpler case. It is the equation satisfied by
the law of a Brownian motion integrated in time. Kolmogorov then wrote
the fundamental solution associated with a Dirac distribution
$\delta_{x_0,v_0}$ initial data:
\begin{align*} 
  G(t,x,v) = 
  \left(\frac{\sqrt 3}{2 \pi t^2}\right)^d
  \exp \left[ - \frac{3 |x-x_0-tv_0 - t(v-v_0)/2|^2}{t^3} -
  \frac{|v-v_0|^2}{4t} \right].
\end{align*}
  
The starting point of H\"ormander's seminal paper \cite{hormander} is
the observation that this fundamental solution shows regularisation in
all variables, even though the diffusion acts only in the velocity
variable. The regularisation in $(t,x)$ is produced by the interaction
between the transport operator $v \cdot \nabla_x$ and the diffusion in
$v$. H\"ormander's paper then proposes general geometric conditions
for this regularisation, called \emph{hypoelliptic}, to hold, based on
commutator estimates. In short, given $X_0$, $X_1$, \dots, $X_n$ a
collection of smooth vector fields on $\R^N$ and the second-order
differential operator $L = - \frac12 \sum_{i=1} ^n X_i^* X_i + X_0$, then
the semigroup $e^{t L}$ is regularising (hypoelliptic) as soon as
the Lie algebra generated by $X_0$, \dots, $X_n$ has dimension $N$
throughout the domain of $L$.

Let us also mention the connexion with the \emph{Malliavin calculus}
in probability, which gives a probabilistic proof to the H\"ormander
theorem in many settings, see \cite{MR536013} as well as the many
subsequent works, for instance
\cite{MR780762,MR783181,MR621660,MR942019}.

\subsection{Extending the DGNM theory to hypoelliptic settings}

The main question of interest here is the extension of the DGNM theory
to hypoelliptic PDEs of divergent type. Hypoelliptic PDEs of second
order $L = - \frac12 \sum_{i=1} ^n X_i^* X_i + X_0$ naturally split
into two classes: the simpler ``type I'' when $X_0=0$ and the operator
is a sum of squares, and the ``type II'' such as the Kolmogorov
equation above, where $X_0 \not = 0$ and the operator combines a
first-order anti-symmetric operator with some partially diffusive
second-order operator. Two main research groups had already been
working on the question.

The extension of the DGNM theory to hypoelliptic operators of ``type
I'' is relatively straightforward. Regarding the ``type II'', Polidoro
and collaborators
\cite{polidoro1994,MR1289901,MR1463798,MR1662349,pr01,pp,MR2368979}
had obtained the Hölder regularity for coefficients with various
continuity assumptions, and had obtained the improvement of
integrability and pointwise bound for measurable coefficients (see
also the isolated result \cite{dfp} for equations in non-divergence
form). Wang and Zhang \cite{wz09,wz11,zhang11} had extended the proof
of Moser for the ``type II'' equations to obtain H\"older regularity,
with technical calculations that did not seem easy to export. Note
also that the use of the DGNM theory in kinetic theory had also been
advocated almost a decade before in the premonitory lecture notes
\cite{Villani-Peccot}.

We present here the work \cite{gimv} (see also the two previous
related preprints \cite{gv,im}) that (1) provides an elementary and
robust proof of the gain of integrability and H\"older regularity in
this ``type II'' hypoelliptic setting, (2) proves the stronger
\emph{Harnack inequality} for these equations (i.e. a quantitative
version of the strong maximum principle).

Let us consider the following kinetic Fokker-Planck equation
\begin{equation}
  \label{eq:main}
  \partial_t f + v \cdot \nabla_x f = \nabla_v \cdot \left(A \nabla_v
    f  \right) + B \cdot \nabla_v f + s, \quad t \in (0,T),\  (x,v) \in \Omega, 
\end{equation}
where $\Omega$ is an open set of $\R^{2d}$, $f = f(t,x,v)$, $B$ and
$s$ are bounded measurable coefficients depending on $(t,x,v)$, and
the $d \times d$ real matrices $A$, $B$ and source term $s$ are
measurable and satisfy
\begin{equation}\label{eq:ellipticity}
 0 < \lambda I \le A \le  \Lambda I, \qquad
 |B| \le  \Lambda, \qquad 
 s \text{ essentially bounded}
\end{equation} 
for two constants $\lambda$, $\Lambda >0$. 
Given $z_0= (t_0,x_0,v_0) \in \R^{2d+1}$, we define the ``cylinder''
$Q_r(z_0)$ centered at $z_0$ of radius $r$ that respects the
invariances of the equation:
\begin{equation}
\label{def:slanted}
Q_r (z_0) := \left\{ (t,x,v) \in \R^{2d+1} \, : \, 
  |x-x_0 - (t-t_0)v_0| < r^3, |v-v_0|< r, 
  t \in \left(t_0-r^2,t_0\right]\right\}.
\end{equation}

The \emph{weak solutions} to equation~\eqref{eq:main} on
$I \times U_x \times U_v$ with $U_x \subset \R^d$ open,
$U_v \subset \R^d$ open, $I = [a,b]$ with $-\infty<a<b \le +\infty$,
are defined as functions
$f \in L^\infty_t(I,L^2_{x,v}(U_x \times U_v))) \cap L^2_{t,x}(I
\times U_x, H^1_v(U_v))$ such that
$\partial_t f + v\cdot \nabla_x f \in L^2_{t,x}(I \times U_x,
H^{-1}_v(U_v))$ and $f$ satisfies the equation~\eqref{eq:main} in the
sense of distributions.

\begin{theorem}[H\"older continuity \cite{gimv}]\label{thm:holder}
  Let $f$ be a weak solution of \eqref{eq:main} in
  $\mathcal Q_0 := Q_{r_0}(z_0)$ and let
  $\mathcal Q_1 := Q_{r_1}(z_0)$ with $r_1 < r_0$. Then $f$ is
  $\alpha$-H\"older continuous with respect to $(x,v,t)$ in
  $\mathcal Q_1$ and
  \[ \| f \|_{C^\alpha(\mathcal Q_1)} \le C \left(\|f \|_{L^2(\mathcal
        Q_0)} + \|s\|_{L^\infty(\mathcal Q_0)}\right)\] for some
  $\alpha \in (0,1)$ and $C>0$ only depending on $d$, $\lambda$,
  $\Lambda$, $r_0$, $r_1$ (plus $z_0$ for $C$).
\end{theorem}
In order to prove such a result, we first prove that $L^2$
sub-solutions are locally bounded; we refer to such a result as an
$L^2-L^\infty$ estimate. We then prove that solutions are H\"older
continuous by proving a lemma which is a hypoelliptic counterpart of
De~Giorgi's ``isoperimetric lemma''.

We moreover prove the Harnack inequality:
\begin{theorem}[Harnack inequality \cite{gimv}]\label{thm:harnack}
  If $f$ is a non-negative weak solution of \eqref{eq:main} in $Q_1(0,0,0)$,
  then
  \begin{equation}\label{eq:harnack}
 \sup_{Q^-} f \le C \left(\inf_{Q^+} f +\|s\|_{L^\infty(Q_1(0,0,0))} \right)
\end{equation}
where $Q^+ := Q_R(0,0,0)$ and $Q^- := Q_R (0,0,-\Delta)$ and $C>1$ and
$R,\Delta \in (0,1)$ are small (in particular
$Q^\pm \subset Q_1(0,0,0)$ and they are disjoint), and universal,
i.e. only depend on dimension and ellipticity constants.
\end{theorem}
Note that using the transformation
$\mathcal{T}_{z_0} (t,x,v)= (t_0+t,x_0+ x + t v_0, v_0+v)$, we get a
Harnack inequality for cylinders centered at an arbitrary point
$z_0=(t_0,x_0,v_0)$. 

Our proof combines the key ideas of De Giorgi and Moser and the
\textbf{velocity averaging method}, which is a special type of
smoothing effect for solutions of the free transport equation
$(\partial_t+v\cdot\nabla_x)f=S$ observed for the first time in
\cite{Agosh,gps} independently, later improved and generalized in
\cite{GLPS,DPL}. This smoothing effect concerns averages of $f$ in the
velocity variable $v$, i.e. expressions of the form
$\int_{\R^d}f(t,x,v)\, \phi(v) \, {\rm d} v$ with, say,
$\phi \in C^\infty_c$. Of course, no smoothing on $f$ itself can be
observed, since the transport operator is hyperbolic and propagates
the singularities. However, when $S$ is of the form
\[
  S=\nabla_v \cdot \left(A(t,x,v)\nabla_vf\right)+s,
\]
where $s$ is a given source term in $L^2$, the smoothing effect of
velocity averaging can be combined with the $H^1$ regularity in the
$v$ variable implied by the energy inequality in order to obtain some
regularity in all directions. A first observation of this type (at the
level of a compactness argument) can be found in \cite{PLLCam};
Bouchut \cite{bouchut} had then obtained quantitative Sobolev
regularity estimates.

Our proof of the $L^2-L^\infty$ gain of integrability follows the
so-called ``De Giorgi-Moser iteration'', see \cite{gimv} where it is
presented in both the equivalent formulations of De Giorgi and of
Moser. We emphasize that, in both approaches, the main ingredient is a
local gain of integrability of non-negative sub-solutions. This latter
is obtained by combining a comparison principle and a fractional
Sobolev regularity estimate following from (1) the velocity averaging
method discussed above and (2) energy estimates. We then prove the
H\"older continuity through a De~Giorgi type argument on the decrease
of oscillation for solutions. We also derive the Harnack inequality by
combining the decrease of oscillation with a result about how positive
lower bounds on non-negative solutions deteriorate with time. It is
worth mentioning here that our ``hypoelliptic isoperimetric argument''
is proved non-constructively, by a contradiction method, whereas the
original isoperimetric argument of De Giorgi is obtained by a
quantitative direct argument. It is an interesting open problem to
obtain such quantitative estimates in the hypoelliptic case.


\section{Conditional regularity of the Landau equation}
\label{sec:lc}

\subsection{Previous works and a conjecture}
The infinite smoothing of solutions to the Landau equation has been
investigated so far in two different settings. On the one hand, it has
been investigated for weak spatially homogeneous solutions
(non-negative in $L^1$ and with finite energy), see \cite{MR2038147}
and the subsequent follow-up papers
\cite{MR2149928,MR2425608,MR2462585,MR2476677,MR2476686,MR1055522,d04,molecules,dvI},
and see also the related entropy dissipation estimates in
\cite{dv,d15}, and see the analytic regularisation of weak spatially
homogeneous solutions for Maxwellian or hard potentials in
\cite{MR2557895}. Furthermore, Silvestre \cite{luislandau} derives an
$L^\infty$ bound (gain of integrability) for spatially homogeneous
solutions in the case of moderately soft potentials without relying on
energy methods. Let us also mention works studying modified Landau
equations~\cite{MR2901061,MR2914961} and the work \cite{MR3599518}
that shows, using barrier arguments, that any weak radial solution to
the Landau-Coulomb equation that belongs to $L^{3/2}$ is automatically
bounded and $C^2$. On the other hand, fewer investigations of the
regularity of spatially heterogeneous solutions have been done,
focusing on the regularisation of classical solutions 
\cite{cdh,MR3191417}.

The general question of conditional regularity hence suggests the
following question in the context of the Landau equation: \medskip

\noindent {\bf Conjecture 1.} Any solutions to the Landau equation
\eqref{MBE}-\eqref{eq:LCop} (with Coulomb interaction $\gamma=-3$) on
$[0,T]$ satisfying \eqref{eq:hydro} is bounded and smooth on $(0,T]$.  \medskip

An important progress has been made by solving a weaker version of
this conjecture when the exponent $\gamma \in (-2,0)$, which
corresponds to \emph{moderately soft potentials}, i.e.
$\gamma + 2s >0$ since here $s=1$. We describe in this section the
different steps and combined efforts of different groups.

\subsection{DGNM theory and local H\"older regularity}
The first step is the work \cite{gimv} already mentioned. A corollary
of Theorem~\ref{thm:holder} is the following:

\begin{theorem}[Local H\"older regularity for the LE \cite{gimv}]
  \label{thm:LE-holder}
  Given any $\gamma \in [-3,1]$, there are universal constants $C>0$,
  $\alpha \in (0,1)$ such that any essentially bounded weak solution
  $f$ of \eqref{MBE}-\eqref{eq:LCop} in
  $(-1,0] \times B_1 \times B_1$ satisfying \eqref{eq:hydro} is
  $\alpha$-H\"older continuous with respect to
  $(t,x,v)\in (-1/2,0] \times B_{1/2}\times B_{1/2}$ and
  \[ \| f \|_{C^\alpha\left(\left(-\frac12,0\right] \times B_{\frac12}
        \times B_{\frac12}\right)} \le C \left(\|f \|_{L^2((-1,0]
        \times B_1 \times B_1)} + \|f\|^{2}_{L^\infty((-1,0] \times
        B_1 \times B_1)} \right).\]
\end{theorem}

Note that this theorem includes the physical case of Coulomb
interactions $\gamma =-3$. The adjective ``universal'' for the
constants refers to their independence from the solution.

\subsection{Maximum principles and pointwise bounds}

This line of research originates in the work of Silvestre both on
the spatially homogeneous Boltzmann (SHBE) and Landau (SHLE) equations
\cite{luis,luis-landau}. These papers build upon the ideas of
``nonlinear maximum principles'' introduced in \cite{MR2989434} in the
case of the Boltzmann collision operator, and upon the so-called
``Aleksandrov-Bakelman-Pucci Maximum Principle'' in the case of the
Landau collision operator, see for instance
\cite{MR1351007,MR2494809}.

The main result of \cite{luis-landau} is:
\begin{theorem}[Pointwise bound for the SHLE]
  Let $\gamma \in [-2,0]$ (moderately soft potentials) and let $f$ be
  a classical non-negative spatially homogeneous solution to the Landau
  equation \eqref{MBE}-\eqref{eq:LCop} on $[0,T] \times \R^d$ for some
  $T >0$, and satisfying the assumptions \eqref{eq:hydro}. Then
  $f \lesssim 1+t^{-3/2}$ with constant depending only on the bounds
  \eqref{eq:hydro}.
\end{theorem}
As noted by the author, this estimate implies quite straightforwardly
existence, uniqueness and infinite regularity for the spatially
homogeneous solution. For the difficult case of very soft potentials
$\gamma \in [-3,2)$, this paper includes a weaker result where the
$L^\infty$ bound depends on a certain weighted Lebesgue norms;
unfortunately it is not yet known how to control such norm along
time. This conceptual barrier, when crossing the ``very soft
potentials threshold'', is reminiscent of the state of the art for the
Cauchy theory in Lebesgue and Sobolev spaces by energy estimates, for
both the spatially homogeneous Boltzmann with long-range interactions
\cite{MR2525118} and Landau equation \cite{MR3375485,MR3158719}.

The pointwise bounds estimates were then extended to the spatially
inhomogeneous case in \cite{cameron2017global}. The main result in
this paper is:
\begin{theorem}[Pointwise bound for the LE]
  Let $\gamma \in (-2,0]$ (moderately soft potentials without the
  limit case) and let $f$ be a bounded non-negative weak solution to
  the Landau equation \eqref{MBE}-\eqref{eq:LCop} on
  $[0,T] \times \R^{2d}$ for some $T >0$, satisfying the assumptions
  \eqref{eq:hydro}. Then $f \lesssim (1+t^{-3/2})(1+|v|)^{-1}$ with
  constant depending only on the bounds \eqref{eq:hydro} (and
  \emph{not} on the $L^\infty$ norm of the solution). Moreover if
  $f_{in}(x,v) \le C_0 e^{-\alpha |v|^2}$, for some $C_0 > 0$ and a
  sufficiently small $\alpha >0$ (depending on $\gamma$ and
  \eqref{eq:hydro}), then $f(t,x,v)\le C_1 e^{-\alpha |v|^2}$ with
  $C_1>0$ depending only on $C_0$, $\gamma$ and the bounds
  \eqref{eq:hydro}.
\end{theorem}

The proof relies on using locally the Harnack inequality in
Theorem~\ref{thm:harnack} adapted to the Landau equation and on
devising a clever change of variable to track how this local estimate
behaves at large velocities. The Gaussian bound is then obtained by
combining existing maximum principle arguments at large velocities
(using that well-constructed Gaussians provide supersolutions at large
$v$) in the spirit of \cite{gamba-panferov-villani-2009}, and the
previous pointwise bound for not-so-large velocities. Finally the
authors remarked that the H\"older regularity estimate of
Theorem~\ref{thm:LE-holder} can be made global using the Gaussian
decay bound.

\subsection{Schauder estimates and higher regularity}

Once the $L^\infty$ norm and the H\"older regularity is under control,
the next step is to obtain higher-order regularity. The classical tool
is the so-called \textbf{Schauder estimates} \cite{MR1545448}. The
purpose of such estimates in general is to show that the solution to
an elliptic or parabolic equation whose coefficients are H\"older
continuous gains two derivatives with respect to the data (source
term, initial data). The gain of the two derivatives is obtained in
H\"older spaces: $C^\delta \to C^{2+\delta}$.

Two works have been obtained independently along this line of
research. The first one \cite{henderson2017c} focuses on the use of
combination of H\"older estimates, maximum principles and Schauder
estimates to obtain conditional infinite regularity for solutions to
the Landau equation with moderately soft potentials
$\gamma \in (-2,0)$. The second one \cite{IM-toy} focuses on the use
of these ingredients to ``break the super-criticality'' of the
nonlinearity for a toy model of the Landau equation. Both these works
develop, in different technical ways, Schauder estimates for this
hypoelliptic equation. The main result in \cite{henderson2017c} is:
\begin{theorem}[Conditional regularity for LE]
  Let $\gamma \in (-2,0)$ (moderately soft potentials without the
  limit case) and let $f$ be a bounded non-negative weak solution to
  the Landau equation \eqref{MBE}-\eqref{eq:LCop} on
  $[0,T] \times \R^{2d}$ for some $T>0$, satisfying the assumptions
  \eqref{eq:hydro} and $f_{in}(x,v) \le C_0 e^{-\alpha |v|^2}$, for
  some $C_0 > 0$ and a sufficiently small $\alpha >0$ (depending on
  $\gamma$ and \eqref{eq:hydro}). Then $f$ is smooth and its
  derivatives have some (possibly weaker) Gaussian decay.
\end{theorem}

Note that: (1) the regularity and decay bounds are uniform in time, as
long as the bounds \eqref{eq:hydro} remain uniformly bounded in time,
(2) further conditional regularity are given in the paper for very
soft potentials $\gamma \in [-3,-2]$ but they require higher
$L^\infty_{t,x}L^1_v(1+|v|)^q$ moments and the constants depend on
time when $\gamma \in [-3,-5/2]$ in dimension $3$, (3) a useful
complementary result is provided by \cite{henderson2017local} where a
local existence is proved in weighted locally uniform Sobolev spaces
and the lower bound on the mass is relaxed by using the regularity to
find a ball where the solution is uniformly positive: the combination
of the two papers provide a conditional existence, uniqueness and
regularity result for moderately soft potentials, conditionally to
upper bounds on the local mass, energy and entropy.  \medskip

The work \cite{IM-toy} considers the toy model: 
\begin{equation}\label{eq:toy}
  \partial f + v \cdot \nabla_x f = \rho[f] \nabla_v \left( \nabla_v f
    + v f \right), \quad \rho[f] := \int_{\R^d} f \dd v, 
\end{equation}
in $x \in \T^d$, $v \in \R^d$, $d \ge 1$. This model preserves the
form of the steady state, the ellipticity in $v$, the non-locality, 
the bilinearity and the mass conservation of the LE. It however
greatly simplifies the underlying hydrodynamic and the maximum
principle structure. Here $H^k(\T^d \times \R^d)$ denotes the standard
$L^2$-based Sobolev space. The main result states (note that solutions
are constructed and not conditional here): 
\begin{theorem}\label{theo:main-toy}
  For all non-negative initial data $f_{in}$ such that
  $f_{in}/\sqrt{\mu} \in H^k(\T^d \times \R^{d})$ with $k > d/2$ and
  satisfying $ C_1 \mu \le f_{in} \le C_2 \mu$ for some $C_1, C_2 >0$,
  there exists a unique global-in-time solution $f$ to \eqref{eq:toy}
  with initial data $f_{in}$ satisfying for all time $t>0$:
  $f(t)/\sqrt{\mu} \in H^k(\T^d \times \R^{d})$ and
  $C_1 \mu \le f \le C_2 \mu$ and $f(t,\cdot,\cdot) \in C^\infty$.
\end{theorem}

Note that the initial regularity exponent $k$ could be relaxed with
more work.  A key step of the proof is the Schauder estimate. It gives
the following additional information on this solution: the
\emph{hypoelliptic H\"older norm} $\mathcal H^\alpha$ (defined below)
of $f/\sqrt \mu$ is uniformly bounded in terms of the $L^2$ norm of
$f_{in}/\sqrt \mu$ for times away from $0$. This norm is defined on a
given open connected set $\cQ$ by
\[ \|g\|_{\mathcal H^\alpha(\cQ)} := \sup_{\cQ} |g| + \sup_{\cQ}
  |(\partial_t + v \cdot \nabla_x )g| + \sup_{\cQ}|D^2_v g| +
  [(\partial_t + v \cdot \nabla_x )g]_{\cC^{0,\alpha}(\cQ)} + [D^2_v
  g]_{\cC^{0,\alpha}(\cQ)}\] where $[\cdot]_{\cC^{0,\alpha}(\cQ)}$ is
a H\"older anisotropic semi-norm, 
i.e. the smallest $C>0$ such that
\[ \forall \, z_0 \in \cQ, \ r >0 \ \mbox{ s.t. } \  Q_r(z_0) \subset \cQ, \quad 
  \|g-g(z_0)\|_{L^\infty(Q_r(z_0))} \le C r^{\alpha} \]
where
\begin{align*}
  Q_r (z_0) & :=\left\{ z : \frac1r (z_0^{-1} \circ z) \in Q_1
              \right\} \\
  & =
  \Big\{ (t,x,v) : t_0 - r^2 < t \le t_0, \, |x-x_0 - (t-t_0)v_0| < r^3, \,
  |v-v_0| < r \Big\}
\end{align*}
and $r z := (r^2 t,r^3 x, r v)$ and
$z_1 \circ z_2 := (t_1+t_2,x_1 +x_2 + t_2 v_1,v_1 + v_2)$.
  



The specific contribution of this work is the study of the Cauchy
problem: the maximum principle provides Gaussian upper and lower
bounds on the solution, and we then provide energy estimates and a
blow-up criterion \emph{\`a la} Beale-Kato-Majda \cite{MR763762}. We
then use the extensions of the DGNM and Schauder theories to control
the $L^\infty_x(H^1_v)$ type norm that governs the blow-up. We prove
H\"older regularity through the method of~\cite{gimv}. We then develop
Schauder estimates following the method of \cite{krylov1996} (see also
\cite{polidoro1994,manfredini,dfp,bb,lunardi1997,radkevich2008,henderson2017c}). New
difficulties arise compared with the parabolic case treated in
\cite{krylov1996} in relation with the hypoelliptic structure and we
develop hypoelliptic commutator estimates directly at the level of
trajectories to solve them. We also borrow some ideas from
hypocoercivity \cite{villani2009memoir} in the proof of the so-called
gradient estimate.

Note that it would be interesting to give a proof of Schauder
estimates for such hypoelliptic equations that is entirely based on
scaling arguments in the spirit \cite{simon} (see also the use of such
scaling arguments in \cite{MR3274562}, in the elliptic-parabolic
case). This might indeed prove useful for generalising such estimates
to the integral Boltzmann collision operator, see the next section.

\section{Conditional regularity of the Boltzmann equation}
\label{sec:nonloc}

\subsection{Previous works and a conjecture}

Short time existence of solutions to \eqref{MBE}-~\eqref{BCO} was
obtained in \cite{amuxy-arma2010} for sufficiently regular initial
data $f_0$. Global existence was obtained in \cite{MR2525118} for
moderately soft potentials in the spatially homogeneous case. In the
next subsections, we present the progresses made so far in the case of
moderately soft potentials: the estimate in $L^\infty$ for $t>0$ was
obtained in \cite{luis}, the local H\"older regularity in
\cite{imbertsilvestre}, and finally the polynomial pointwise decay
estimates in \cite{IMS}. The bootstrap mechanism to obtain higher
regularity through Schauder estimates remains however unsolved at now.
 
Let us briefly review the existing results about regularisation. The
very first mathematical observation that long-range interactions are
associated with fractional ellipticity in the kinetic variable goes
back to Desvillettes \cite{MR1324404} in the mid 1990s. In
\cite{amuxy-arma2010}, the authors prove that if the solution $f$ has
five derivatives in $L^2$, with respect to all variables $t$, $x$ and
$v$, weighted by $(1+|v|)^q$ for arbitrarily large powers $q$, and in
addition the mass density is bounded below, then the solution $f$ is
$C^\infty$. It is not known however whether these hypotheses are
implied by \eqref{eq:hydro}. Note also the previous partial results
\cite{MR2038147,MR2149928,MR2425608,MR2462585,MR2476677,MR2476686,MR2885564}
in the spatially homogeneous case and with less assumptions on the
initial data, and the work \cite{MR2820356} in the spatially
inhomogeneous case but with much stronger a priori assumptions.

Note that, drawing inspiration from the case of the Landau equation,
in order for the iterative gain of regularity in
\cite{henderson2017c,IM-toy} to work, it is necessary to start with a
solution that decays, as $|v| \to \infty$, faster than any algebraic
power rate $|v|^{-q}$. We expect the same general principle to apply
to the Boltzmann equation, even if the appropriate Schauder type
estimates for kinetic integro-differential equations to carry out an
iterative gain in regularity are not yet available.


The question of conditional regularity suggests the following
conjecture in the context of the Boltzmann equation with long-range
interactions: \medskip

\noindent {\bf Conjecture 2.} Any solutions to the Boltzmann equation
\eqref{MBE}-\eqref{BCO} with long-range interactions
($\gamma \in (-3,1]$, $s \in (0,1)$, $\gamma + 2s \in (-1,1)$) on
$[0,T]$ satisfying \eqref{eq:hydro} is bounded and smooth on $(0,T]$.
\medskip

The rest of this section is devoted to describing the partial
progresses made in the case of, again, moderately soft potentials
$\gamma +2s >0$.

\subsection{Maximum principle and pointwise $L^\infty$ bound}

This first breakthrough is due to Silvestre \cite{luis}. This article
draws inspiration from his own previous works on non-local operators
and from the ``nonlinear maximum principle'' of Constantin and Vicol
\cite{MR2989434}. It is based on a maximum principle argument for a
barrier supersolution that is constant in $x,v$ and blowing-up as
$t \to 0^+$; it uses the decomposition of the collision operator and
``cancellation lemma'' going back to \cite{advw}, the identification
of a cone of direction for $(v'-v)$ is order to obtain lower bounds on
the $f$-dependent kernel of the elliptic part of the operator, and
finally some Chebycheff inequality and nonlinear lower bound on the
collision integral. The main result is:
\begin{theorem}[Pointwise bound for the BE \cite{luis}]
  Let $\gamma \in [-2,1]$, $s \in (0,1)$ with $\gamma + 2s >0$
  (moderately soft potentials). Let $f$ be a classical non-negative
  solution to the Boltzmann equation \eqref{MBE}-\eqref{eq:LCop} on
  $[0,T] \times \T^d \times \R^d$ for some $T>0$, satisfying the
  assumptions \eqref{eq:hydro}. Then $f \le C( 1+t^{-\beta})$ with
  $C>0$ and $\beta >0$ depending only on $\gamma$, $s$ and the bounds
  \eqref{eq:hydro}.
\end{theorem}

Note that the paper also includes further results in the case of very
soft potentials but conditionally to additional estimates of the form
$L^\infty_{t,x} L^p_v(1+|v|^q)$ for some $p>1$, $q>0$; it is not known
at present how to deduce the latter estimates from the hydrodynamic
bounds \eqref{eq:hydro}. 

\subsection{Weak Harnack inequality and local H\"older regularity}

The second breakthrough is the paper \cite{imbertsilvestre}. In
comparison to the Landau equation, the Boltzmann equation has a more
complicated integral structure, that shares similarity with ``fully
nonlinear'' fractional elliptic operators. The main result proved is:
\begin{theorem}[Local H\"older regularity for the BE \cite{imbertsilvestre}]
  Given any $\gamma \in (-3,1]$ and $s \in (0,1)$ with
  $\gamma + 2s >0$, there are universal constants $C>0$,
  $\alpha \in (0,1)$ such that any essentially bounded non-negative
  weak solution $f$ to \eqref{MBE}-\eqref{BCO} in
  $(-1,0] \times B_1 \times \R^3$ satisfying \eqref{eq:hydro} is
  $\alpha$-H\"older continuous with respect to
  $(t,x,v)\in (-1/2,0] \times B_{1/2}\times B_{1/2}$, where
  $C, \alpha$ only depend on the $L^\infty$ bound of $f$ and the
  bounds \eqref{eq:hydro}.
\end{theorem}

The proof goes in two steps. The first step is a local
$L^2 \to L^\infty$ gain of integrability, following the approach of De
Giorgi and Moser as reformulated in a kinetic context in \cite{pp} and
\cite{gimv}. It requires further technical work to formulate the De
Giorgi iteration for such integro-differential equations with
degenerate kernels (see also the related works \cite{K,K1,K2}). The
regularity mechanism at the core of the averaging velocity method is
used; it is however presented differently than in most papers on this
topic, by relying on explicit calculations on the fundamental solution
of the fractional Kolmogorov equation. In the second step of the
proof, the authors establish a weak Harnack inequality, i.e. the
control from above of local $L^\epsilon_{t,x,v}$ averages with
$\epsilon >0$ small by a local infimum multiplied by a universal
constant. This inequality is sufficient to deduce the H\"older
regularity. Two different strategies are used depending on whether
$s \in (0,1/2)$ or $s \in [1/2,1)$. In the first case, they construct
a barrier function to propagate lower bounds as in the method by
Krylov and Safonov for equations in nondivergence form. In the second
case, they use a variant of the isometric argument of De Giorgi proved
by compactness as in \cite{gimv}. Again the regularity of velocity
averages plays a crucial role but is exploited by direct calculation
on the fundamental solution of the fractional Kolmogorov equation.

\subsection{Maximum principle and decay at large velocities}

Finally in the paper \cite{IMS}, the nonlinear maximum principle
argument of \cite{luis} is refined to obtain ``pointwise counterpart''
of velocity moments. 
The main result established in this paper is:
\begin{theorem}[Decay at large velocities for the BE \cite{IMS}]
  Given any $\gamma \in (-3,1]$ and $s \in (0,1)$ with
  $\gamma + 2s \ge 0$, there are universal constants $C>0$,
  $\alpha \in (0,1)$ such that for any classical non-negative solution
  $f$ to \eqref{MBE}-\eqref{BCO} in $[0,T] \times \T^3 \times \R^3$
  satisfying \eqref{eq:hydro}, it holds for any $q >0$: (i) 
  if $f_{in}\lesssim (1+|v|)^{-q}$ then $f (t,x,v) \le C (1+|v|)^{-q}$
  for all $t>0$, (ii) assuming furthermore that $\gamma>0$, then
  $f(t,x,v) \le C' (1 + t^{-\beta}) (1 + |v|)^{-q}$ for all $t>0$. All
  the constants depend on $\gamma$, $s$, $q$ and the bounds
  \eqref{eq:hydro}.
\end{theorem}

The study of large velocity decay in weighted $L^1$ spaces, known as
the study of \emph{moments}, is an old and important question in
kinetic equations.  The study of moments was initiated, for spatially
homogeneous solutions, in \cite{MR0075725} for Maxwellian potentials
($\gamma=0$). In the case of hard potentials ($\gamma>0$), Povzner
identities \cite{MR0142362,MR684411,MR1461113,MR1478067} play an
important role.  For instance, Elmroth \cite{MR684411} used them to
prove that if moments are initially bounded, then they remain bounded
for all times.  Desvillettes \cite{MR1233644} then proved that only
one moment of order $q>2$ is necessary for the same conclusion to hold
true. It is shown in \cite{MR1461113,MR1697562} that even the
condition on one moment of order $q>2$ can be dispensed with, in both
(homogeneous) cutoff and non-cutoff case. These moment estimates were
used by Bobylev \cite{MR1478067} in order to derive (integral)
Gaussian tail estimates. In the case of soft potentials, Desvillettes
\cite{MR1233644} proved for $\gamma \in (-1,0)$ that initially bounded
moments grow at most linearly with time and it is explained in
\cite{review} that the method applies to $\gamma \in [-2,0)$.  The
case of measure-valued solutions is considered in \cite{MR2871802}.

However the extension of these integral moments estimates to spatially
inhomogeneous solutions is a hard and unclear question at the
moment. The only result available is \cite[Lemma~5.9 \& 5.11]{GMM}
which proves the propagation and appearance of certain exponential
moments for the spatially inhomogeneous Boltzmann equation for hard
spheres (or hard potentials with cutoff), however in a space of the
form $W^{3,1}_x L^1(1+|v|^q)$. Another line of research opened by
\cite{gamba-panferov-villani-2009} consists in establishing
exponential Gaussian pointwise decay by maximum principle arguments
(see also
\cite{bobylev2017upper,alonso2017exponentially,gamba2017pointwise}). But
these works assumes exponential integral moments whose propagation in
time is not known, therefore it is not clear how to use them in this
context.

We finally recall that the last part of the research program, the
Schauder estimates, is missing for the Boltzmann equation with
moderately soft potentials, and is an interesting open question for
future research.

\bigskip

\noindent {\bf Acknowledgement.} This review article greatly benefited
from the many discussions of the author with Cyril Imbert and Luis
Silvestre, in particular concerning the works on the conditional
regularity initiated by Luis Silvestre that form a key part of this
review, and concerning their important joint work on the weak Harnack
inequality for the Boltzmann equation. Moreover the author gratefully
acknowledges the many crucial discussions with Alexis Vasseur about
his insights on the De Giorgi method, and with Fran\c{c}ois Golse
about averaging lemma and homogeneisation. The author acknowledges
partial funding by the ERC grants MATKIT 2011-2016 and MAFRAN
2017-2022.

\bibliographystyle{acm}
\bibliography{icm}

\bigskip

\signcm 

\end{document}